\documentclass[11pt]{article}
\usepackage{epsfig}
\usepackage{amssymb}
\usepackage[]{times}

\makeatletter
\@addtoreset{equation}{section}

\makeatother
\newcommand{\be}{\begin{equation}}
\newcommand{\ee}{\end{equation}}
\newcommand{\bea}{\begin{eqnarray}}
\newcommand{\eea}{\end{eqnarray}}

\begin{document}

We present English translation of the classical article of Hermann Amadeus 
Schwarz (1843--1921)\\"{\em Proof of the theorem that a surface area of a ball 
is smaller than of any other body of the same volume}" which was published in 
1884, in {\sf Proceedings of the K\"{o}nigliche Gesellschaft der Wissenschaften 
and the Georg-Augusts-Universit\"{a}t, G\"{o}ttingen}.

We preserved the author notations throughout the text and tried to follow his 
grammar construction of the sentences. One editorial comment in the footnote at 
the page 4 is related to the (possible) misprint in the original German text. 
The other three comments in the footnotes at the pages 6, 7, 8 are given to 
warn about ambiguities in definition of different notions and entities.

\vspace{1cm}
\begin{verbatim}
Simon Ulka,
International School,
Technion - Israel Institute of Technology,
32000 Haifa, Israel
e-mail: simon@ulka.de

Leonid G. Fel,
Department of Civil Engineering,
Technion  - Israel Institute of Technology,
32000 Haifa, Israel  
e-mail: lfel@technion.ac.il

Boris Y. Rubinstein,
Stowers Institute for Medical Research,
1000 E 50th St,  
Kansas City, MO 64110, USA
e-mail: bru@stowers.org
\end{verbatim}

\newpage
\centerline{\LARGE\bf Proof of the theorem that a surface area of a ball}
\vskip0.1 cm
\centerline{\LARGE\bf is smaller than}
\vskip0.1 cm
\centerline{\LARGE\bf of any other body of the same volume}
\vskip0.2 cm
\centerline{\LARGE\sf by H. A. Schwartz}
\vskip0.4 cm
\centerline{\sf\large Proceedings of the K\"{o}nigliche Gesellschaft der
Wissenschaften}
\centerline{\sf\large and the Georg-Augusts-Universit\"{a}t, G\"{o}ttingen,
1884, pages 1-13}
\vskip0.5 cm

\centerline{\large\sl Translated by S. Ulka}
\vskip0.1 cm

\centerline{\large\sl Edited by L. Fel and B. Rubinstein}
\vskip0.5 cm

In order to prove the theorem that {\it a surface area of a ball is smaller
than that of any other body of the same volume}, several different approaches
were used, which are mainly based on the condition that {\it among all bodies
of same volume, there exists one that has a minimal surface area}. As long as
the condition, is not proved none of the aforementioned methods are valid to be
used in order to prove the main theorem.

Trying to prove the above mentioned {\it theorem} for bodies which surface is
formed by a finite number of finite analytical surface pieces, I was led to a
method which does not seem to be exposed to the objection of missing rigor. The
proof presented in this manuscript, is based on the repetitive use of an method 
that has been used by Mr. {\it Weierstrass} in his lectures on calculus of 
variation. I'm indebted to one of his oral presentations for the knowledge 
about this method.

\centerline{\S 1}

\vspace{.1cm}
Let ${\mathfrak U}\;$ be a non-spherical body which surface ${\mathfrak B}$ is
formed by a finite number of pieces of analytical surfaces assumed to be free
of singularities.

The points of the surface ${\mathfrak B}$ shall be related to a right-angular
coordinate system, which is chosen such that no part of ${\mathfrak B}$ is
parallel to the $yz$-plane of the coordinate system. Let $x_0$ be the smallest
and $x_1$ be the largest of all accessible values of the coordinate $x$.

Let an arbitrary point $P$ that belongs to ${\mathfrak B}$ and does not belong
to any of its {\it edges} has a right-angular coordinates $x,y,z$. Construct a
normal to ${\mathfrak B}$ at $P$ and fix its positive direction such that at
point $P$ its direction inclines from the outside to the inside w.r.t. the body
${\mathfrak U}\;$. The angle between the positive direction of the normal and
the positive direction of the $x$-axis shall be denoted by $\xi$.

Through this point $P$ construct a plane ${\mathfrak E}_x$ which is parallel to
the $yz$-plane. This plane has generally one or more curves ${\mathfrak C}$ in
common with the surface ${\mathfrak B}$. A totality ${\mathfrak C}_x$ of these
curves shall be considered as one curve, while $dy$, $dz$ stand for coordinates
of an element of ${\mathfrak C}_x$ coming from point $P$ with length $ds$. One
can determine that the initial points for measuring the arc length $s$ on the
curves ${\mathfrak C}_x$ ($x_0<x<x_1$) form one or more analytical curves that
belong to the surface ${\mathfrak B}$. 

The values of $x$ and $s$ shall be chosen as independent variables to define 
the coordinates $x,y,z$ of any point on the surface ${\mathfrak B}$.

The abovementioned assumptions imply that it is always possible to cut the
surface ${\mathfrak B}$ by a finite number of planes parallel to the $yz$-plane
into a finite number of either bowl-shaped or ring-shaped partial surfaces,
such that for each of them the coordinates $y$ and $z$ of any point are unique
and generally continuous functions of two variables $x$ and $s$. If the curve
${\mathfrak C}_x$ is made of several separate closed curves, then one has to
assign a specific order as well as starting points for the measurement of the
arc length. The value of $s$ shall be defined in such a way that if $U_1(x),
U_2(x),\ldots,U_n(x)$ is the {\it length} of the first, second, \ldots, $n$-th
closed parts of the curve ${\mathfrak C}_x$, then $s$ takes all the values from
$0$ up $U_1(x)$ for the first part, $U_1(x)$ to $U_1(x)+U_2(x)$ for the second,
\ldots, $U_1(x)+U_2(x)+\ldots+U_{n-1}(x)$ to $U(x)$ for the last one, where
$U(x)$ denotes the total length of all parts of ${\mathfrak C}_x$.

Therefore the following equations hold:
$$
\frac{\partial x}{\partial s}=0,\qquad
\left(\frac{\partial y}{\partial s}\right)^2+
\left(\frac{\partial z}{\partial s}\right)^2=1.
$$

Assume that a growth of $s$ along each segment of the curve ${\mathfrak C}_x$
was chosen such that
$$
A =\frac{\partial y}{\partial x}\frac{\partial z}{\partial s}-
\frac{\partial z}{\partial x}\frac{\partial y}{\partial s},\qquad
B = -\frac{\partial z}{\partial s},\qquad
C = \frac{\partial y}{\partial s},
$$
where $A,B,C$ are coordinates of a line segment which direction is the same as
the direction of the {\it normal} to the surface ${\mathfrak B}$ for every 
point $P$ that does not belong to any {\it edge} of the surface ${\mathfrak B}$.

Under this condition, the integral along the curve ${\mathfrak C}_x$
$$
\int_{0}^{U(x)} \frac{1}{2}\left(y\frac{\partial z}{\partial s}-
z\frac{\partial y}{\partial s}\right)ds = Q(x)
$$
is the area of the surface that consists of one or more pieces, which include
all points that are inside of body ${\mathfrak U}$ and the plane ${\mathfrak E}
_x$ and none else.

The volume $V$ of the body ${\mathfrak U\;}$ is given by following equation
$$
\int_{x_0}^{x_1} Q(x) dx = V.
$$
The size of an element $dS$ on the surface ${\mathfrak B}$ and the size $dT$
of an orthogonal projection of $dS$ onto the $yz$-plane of the coordinate
system are given by equations
\footnote{In the original German text instead of $\sqrt{1+A^2}$ there is
erroneously written $\sqrt{1+A}$.}
$$
dS = \sqrt{1+A^2}\;dxds=\frac{1}{\sin\xi}dxds,\qquad
dT = \cos\xi\;dS =A dxds= \cot\xi\;dxds.
$$
From the geometric meaning of $Q(x)$ and the integral along the curve
${\mathfrak C}_x$
$$
\int_{s=0}^{s=U(x)}\cot \xi\; dxds =dx\int_{0}^{U(x)}\cot \xi\;
ds= dx\int_{0}^{U(x)}A ds
$$
we obtain equation
$$
dQ(x) = Q'(x)dx = dx\int_{0}^{U(x)}A ds.
$$

The area $S$ of the surface ${\mathfrak B}$ and the volume $V$ of the body
${\mathfrak U}\;$ therefore read
\bea
&&S=\int_{x_0}^{x_1}dx\int_{0}^{U(x)}\sqrt{1+\left(
\frac{\partial y}{\partial x}\frac{\partial z}{\partial s}-
\frac{\partial z}{\partial x}\frac{\partial y}{\partial s}\right)^2} \;ds=
\int_{x_0}^{x_1}dx\int_{0}^{U(x)}\sqrt{1+A^2}\;ds,\nonumber\\
&&V = \int_{x_0}^{x_1}dx\int_{0}^{U(x)}\frac{1}{2}
\left(y\frac{\partial z}{\partial s}-z\frac{\partial y}{\partial s}\right)ds=
\int_{x_0}^{x_1} Q(x) dx.\nonumber
\eea
\vspace{.1cm}

\centerline{\S 2}

\vspace{.2cm}
There exists a possibility that for one or more values of $x$ the entity
$$
\frac{\partial y}{\partial x}\frac{\partial z}{\partial s}-
\frac{\partial z}{\partial x}\frac{\partial y}{\partial s}=A
$$
reaches a value $Q'(x)/U(x)$ which is independent of $s$.

Under this condition there exists $x$ such that
$$
\int_{0}^{U(x)}\sqrt{1+A^2}\;ds = \sqrt{U^2(x)+Q'^2(x)}.
$$
But in any other case it holds that
$$
\int_{0}^{U(x)}\sqrt{1+A^2}\;ds > \sqrt{U^2(x)+Q'^2(x)}.
$$
In order to prove this statement, set
$$
\int_{0}^{s}Ads = t,\qquad
\int_{0}^{s}\sqrt{1+A^2}ds = \int_{0}^{s}\sqrt{ds^2+dt^2}
$$
and determine an angle $\omega$ by equation
$$
\cos\omega = \frac{sds+tdt}{\sqrt{s^2+t^2}\sqrt{ds^2+dt^2}}
$$
resulting in formula
$$
d\left(\int_{0}^{s}\sqrt{ds^2+dt^2} -\sqrt{s^2+t^2}\right) =
(1-\cos\omega)\sqrt{ds^2+dt^2},
$$
which geometric interpretation is self-explanatory and therefore does not
require any further remarks.

Since $(1-\cos\omega)$ is non-negative and can vanish everywhere in the
interval $0<s<U(x)$, only when $A$ is independent of $s$, integrating both
sides of the equation (given above) between $s=0$ and $s=U(x)$ one obtains
$$
\mbox{(I.a)}
\quad\quad\quad\int_{0}^{U(x)}\sqrt{1+A^2}ds \ge \sqrt{U^2(x)+Q'^2(x)}.
$$
Here the {\it equality} can be reached only when the values of $x,\;A,$ and,
therefore, $\xi$ are independent of $s$.

(I.a) implies
$$
\mbox{(I.b)}
\quad\quad\quad\quad\quad S \ge \int_{x_0}^{x_1} \sqrt{U^2(x)+Q'^2(x)}\; dx.
$$
The {\it equality} can  be reached only  when in the entire interval
$x_0<x<x_1$ it holds that $A = \cot \xi$ is a function of $x$ alone.
\vspace{.2cm}

\centerline{\S 3}

\vspace{.1cm}
It may occur that for one or more values of $x$ the closed curve ${\mathfrak C}
_x$ is formed by a single {\it circular arc} of radius $r$. Under this
condition we obtain for values of $x$:
$$
r^2\pi=Q(x),\qquad U^2(x)=(2r\pi)^2 = 4Q(x)\pi,\qquad\sqrt{U^2(x)+Q'^2(x)}=
\sqrt{4Q(x)\pi+Q'^2(x)}.
$$
In any other case though it holds that
$$
U^2(x) \ge 4Q(x)\pi.
$$

To prove this assertion, in other words, to show that {\it the perimeter of a
circular surface is smaller than that of any other planar figure of the same
area}, we can proceed as follows.

In order to handle the case when the curve ${\mathfrak C}_x$ consists of several
closed curves similarily as in the case when ${\mathfrak C}_x$ is given by a
single closed curve, one can assume that the integral taken along a given
{\it closed} curve, that lies inside the plane ${\mathfrak E}_x$,
$$
\int\frac{1}{2}\left(y dz-z dy\right)
$$
does not change its value, when the curve is translated in the plane without
changing its shape. Since $ds=\sqrt{dy^2+dz^2}$ does not change during this
operation, it can be used to translate all closed curves produced by
intersection of the plane ${\mathfrak E}_x$ with surface ${\mathfrak B}$,
which produce ${\mathfrak C}_x$. We translate these curves in such a manner
that the point on each curve, where $s$ has its largest or smallest value, is
the same as the point $O_x$ of intersection of the plane ${\mathfrak E}_x$ with
the $x$-axis.

By moving the pieces of the curve ${\mathfrak C}_x$ and combining them into a
single polygon chain -- such that along this chain $s$ grows monotonically,
taking all values in the interval $0 < s < U(x)$ -- one produces a single
closed curve, which shall be called $\bar{\mathfrak C}_x$. The length of this
curve is $U(x)$.

If we denote by $y$ and $z$ the second and the third coordinates of any point
$P$ on $\bar{\mathfrak C}_x$ that is determined by the value of $s$, then these
variables are unique and {\it continuous} functions of $s$ in the interval $0<s
< U(x)$ and equal to zero at the ends of this interval.

The integral along the curve $\bar{\mathfrak C}_x$
$$
\int_{0}^{U(x)}\frac{1}{2}
\left(y\frac{\partial z}{\partial s}-z\frac{\partial y}{\partial s}\right)ds
$$
has a value of $Q(x)$.

We shall now introduce notations
$$
\rho = \sqrt{y^2+z^2},\quad {\mathfrak F} = \int_{0}^{s}\frac{1}{2}
\left(y\frac{\partial z}{\partial s}-z\frac{\partial y}{\partial s}\right)ds,
$$
and consider a {\it circular segment} with the chord belonging to the plane
${\mathfrak E}_x$ along the line segment $O_xP$, which is constructed in such
manner that along the segment arc in the direction from $O_x$ to $P$ the
integral
$$
\int_{0}^{s}\left(y dz-z dy\right)
$$
has the value ${\mathfrak F}$.

Denote,

$r$, positive or negative (according to the sign of ${\mathfrak F}$) value of
the radius of the 
circular segment\footnote{In beginning of \S 3 the variable $r$ 
is defined by equality $r^2\pi=Q(x),$ as a radius of circular arc
which forms the closed curve ${\mathfrak C}_x$.};

$\phi$, positive or negative value of a half of the central angle of
the circular segment measured in radians;

${\mathfrak L}$, length of the arc. The following equations hold:
$$
{\mathfrak F} = r^2(\phi-\sin\phi\cos\phi),\qquad\rho=2r\sin\phi,\qquad
{\mathfrak L} = 2r \phi.
$$

To simplify the analysis, we can assume, without loss of generality, that both
$\rho$ and ${\mathfrak F}$ will never equal to $0$ {\it simultaneously} for
values of $s$ {\it inside} the interval $0 < s < U(x)$. In the case when both
$\rho$ and ${\mathfrak F}$ equal to $0$ at $s_0$ in the interval $0<s<U(x),$
but are not $0$ at the same time for any value of $s$ in the interval $s_0<s<
U(x)$, we can limit the analysis of the interval $0 < s < U(x)$ to the interval
$s_0 < s < U(x)$. Then it holds that
$$
U(x) > \sqrt{4Q(x)\pi}
$$
a fortiori.

To determine the value of $\phi$ in the interval $-\pi\leq\phi\leq\pi$, we
obtain a transcendent equation
$$
f(\phi) = \frac{\phi-\sin\phi\cos\phi}{2\sin^2\phi}=
\frac{2{\mathfrak F}}{\rho^2}.
$$
Since the derivative
$$
f'(\phi) = \frac{\tan\phi-\phi}{\tan\phi\sin^2\phi}
$$
is neither negative nor $0$ for all values of $\phi$ in the above interval,
while the function $f(\phi)$ takes all values from $-\infty$ to $+\infty$ in
that interval, then for every value of $s$ in the given interval the equation
$$
f(\phi) = \frac{2{\mathfrak F}}{\rho^2}
$$
has only one root which value continuously changes with the value of $s$.

The variable $r$ that has the same sign as $\phi$ and ${\mathfrak F}$ is thus
defined completely. The value of $1/r$ changes continuously with $s$.

The curve $\bar{\mathfrak C}_x$ intersects with the line segment $O_xP$ as
well as with the arc of the previously constructed circular segment at point
$P$. If we denote by $\omega$ an angle\footnote{In \S 2 $\omega$ is defined by 
equation $\cos\omega=(sds+tdt)/(\sqrt{s^2+t^2}\sqrt{ds^2+dt^2})$.} formed by 
the tangents to the curve $\bar{\mathfrak C}_x$ and to the circular arc that 
connects $O_x$ and $P$ at the point $P$ then following equations hold:
\bea
\cos\omega & = &\frac{1}{\rho}\left[\left(y\frac{\partial y}{\partial s}+
z\frac{\partial z}{\partial s}\right)\cos\phi+
\left(y\frac{\partial z}{\partial s}-z\frac{\partial y}{\partial s}
\right)\sin\phi\right],\nonumber\\
\sin\omega & = &\frac{1}{\rho}\left[\left(y\frac{\partial y}{\partial s}+
z\frac{\partial z}{\partial s}\right)\sin\phi-
\left(y\frac{\partial z}{\partial s}-z\frac{\partial y}{\partial s}
\right)\cos\phi\right],\nonumber\\
d{\mathfrak F}&=&\frac{1}{2}\left(y\frac{\partial z}{\partial s}-
z\frac{\partial y}{\partial s}\right)ds,\quad
d\rho=\frac{1}{\rho}\left(y\frac{\partial y}{\partial s}+
z\frac{\partial z}{\partial s}\right)ds,\nonumber \\
d{\mathfrak L}&=& \cos\omega \;ds,\quad
\frac{\tan\phi-\phi}{\tan\phi\sin^2\phi}d\left(1/r \right) =
-\frac{2\sin\omega}{\rho^2} ds.\nonumber
\eea
We therefore obtain an equation:
$$
d(s-{\mathfrak L})=(1-\cos\omega)ds.
$$
Since $s$ takes all values from $0$ to $U(x)$ while continuously growing and
the term $1-\cos\omega$ is nonnegative and only under a certain condition it
{\it everywhere} equals $0$, an equation holds
$$
s-{\mathfrak L} = \int_{0}^{s}(1-\cos\omega)ds,
$$
from which we deduce that the {\it length} $s=U(x)$ of the curve $\bar
{\mathfrak C}_x,$ is generally {\it larger} than ${\mathfrak L}$ that
corresponds to the value $s=U(x)$.

For $s=U(x)$, $\rho$ turns to $0$, ${\mathfrak F}$ reaches $Q(x)$, $\phi$
becomes $\pi$ because $Q(x)$ is positive and the circular segment is replaced
by a circular surface with radius $\sqrt{Q(x)/\pi}$, and ${\mathfrak L}$
equals  $\sqrt{4Q(x)\pi}$.

From the previous analysis we find following relationship:
$$
\mbox{(II.a)}\hspace{2cm} U(x)\ge\sqrt{4Q(x)\pi}.
$$

It is important to note that, in accordance with the statement made above, the
{\it equality} takes place if and only if in the interval $0<s<U(x)$ the
quantity $r$ has a value of $\sqrt{Q(x)/\pi},$ independent of $s$. It means
that the curve ${\mathfrak C}_x$ which belongs to the plane ${\mathfrak E}_x$
is a {\it circular arc} with radius $\sqrt{Q(x)/\pi}$. Only under this condition
the value of the terms $1-\cos\omega$ and $\sin\omega$ is $0$ in the whole
interval $0<s<U(x)$.

Combining formula (I.b) with (II.a) we find:
$$
\mbox{(II.b)}\hspace{2cm}S\ge\int_{x_0}^{x_1}\sqrt{4Q(x)\pi+Q'^2(x)}\;dx.
$$
It should be noted that the {\it equality} is reached only when the body
${\mathfrak U}$ is a {\it body of rotation} with rotation axis parallel to the
$x$-axis of the coordinate system.
\vspace{.1cm}

\centerline{\S 4}

\vspace{.1cm}
Consider a {\em body of rotation} ${\mathfrak D}$, which rotation axis coincides
with the $x$-axis of the coordinate system, with a radius $r$ of the parallel
circle that lies in the plane ${\mathfrak E}_x$ given by equations
$$
r^2\pi = Q(x),\qquad r \ge 0,\qquad x_0 \le x \le x_1.
$$
Then it appears that
$$
2r\pi dr=Q'(x) dx,\qquad 2r\pi\sqrt{dx^2+dr^2}=\sqrt{4Q(x)\pi+Q'^2(x)}\; dx.
$$
The surface area of the body of rotation ${\mathfrak D}$ is given by integral
$$
\int_{x_0}^{x_1}\sqrt{4Q(x)\pi+Q'^2(x)}\; dx.
$$
Introduce the notation
\footnote{In \S 1 ${\mathfrak B}$ denotes a surface of a non-spherical body
${\mathfrak U}\;$.}
$$
{\mathfrak B} =\int_{x_0}^{x_1}Q(x)\;dx
$$
and consider a {\em spherical segment} that lies on the negative side of the
plane ${\mathfrak E}_x,$ which volume is equal to ${\mathfrak B}$. The plane
boundary has an area $r^2\pi = Q(x)$ and coincides with the parallel circle of
${\mathfrak D}$. The area of the curved boundary of the spherical segment, the
so called {\em spherical cap}, shall be called ${\mathfrak H},$ its radius is
$R$, and one quarter of the central angle in radians is $\psi$. The angle that
the tangential plane to the spherical surface makes with the tangential plane to
the body of rotation in the point of the common parallel circle will be denoted
by $\omega$ and the length $dl$ of the meridian segment of the body of rotation
is equal to $\sqrt{dx^2 + dr^2}$.

With these assumptions, following equations hold
\bea
{\mathfrak B} &=&\frac{4}{3}R^3\pi\sin^4\psi(\sin^2\psi+3\cos^2\psi),\nonumber\\
r&=&2R\sin\psi\cos\psi,\nonumber\\
{\mathfrak H}&=&4R^2\pi\sin^2\psi,\nonumber\\
d{\mathfrak B}&=&r^2\pi dx,\qquad d{\mathfrak H}=2r\pi\cos\omega dl,\nonumber\\
\cos\omega&=&\cos 2\psi\frac{dr}{dl}+\sin 2\psi\frac{dx}{dl},\nonumber\\
\sin\omega&=&\cos 2\psi\frac{dx}{dl}-\sin 2\psi\frac{dr}{dl}.\nonumber\\
(1+\tan^2\omega)^2d\left(\frac{1}{R}\right)&=&\frac{4\sin\psi}{r^2}dl,\nonumber
\eea
The value of $\psi$ can be determined from the equation
$$
f(\psi) = \tan\psi + \frac{1}{3}\tan^3\psi = \frac{2B}{r^3\pi},
$$
under the condition that
$$
0 \le \psi \le \pi/2.
$$
Since the derivative of the function $f(\psi)$
$$
f'(\psi) = (1+\tan^2\psi)^2
$$
is always positive, the function $f(\psi)$ takes all values from $0$ to
$\infty$ while growing continuously, whereas $\psi$ takes all values from $0$
to $\pi/2$ growing continuously as well. For this reason the equation
$$
f(\psi)=\frac{2{\mathfrak B}}{r^3\pi}
$$
in the given interval has one and only one root, which value continuously
changes with $x$.

Since the value of $\psi$ is unique, we can say the same about $R$ due to the
equation $r=2R\sin\psi\cos\psi$. From the equation
$$
d\left[\int_{x_0}^{x}2r\pi\sqrt{1+(dr/dx)^2} dx - {\mathfrak H}\right]
=2r\pi(1-\cos\omega)dl,
$$
we obtain
$$
\int_{x_0}^{x}\sqrt{4Q(x)\pi+Q'^2(x)} dx \ge {\mathfrak H}.
$$
If now we set $x=x_1$, then $r$ turns to $0$, $\psi$ reaches $\pi/2$,
${\mathfrak B}$ becomes $V$, $R$ results in the value of $\sqrt[3]{3V/4\pi}$
and ${\mathfrak H}$ becomes $\sqrt[3]{36V^2\pi}$.

Furthermore we obtain
$$
\mbox{(III.a)}\quad\quad\quad\quad\quad
\int_{x_0}^{x_1}\sqrt{4Q(x)\pi+Q'^2(x)} dx \ge \sqrt[3]{36V^2\pi}.
$$
The {\it equality} is reached only if in the interval $x_0<x<x_1$, $\cos\omega=
1$, $\sin\omega=0$ and, therefore, $R$ has a value independent of $x$. In this
case though it follows from the above equations that the body of rotation
${\mathfrak D}$ is a {\it ball}.

If that happens, then the body ${\mathfrak U}$ cannot be a body of rotation
which axis of rotation is parallel to the $x$-axis of the coordinate system;
${\mathfrak U}$ would have to be a ball as well ( this possibility was already
eliminated in the beginning).

Combining formulas (II.b) and (III.a), it appears that
$$
S \ge \int_{x_0}^{x_1}\sqrt{4Q(x)\pi+Q'^2(x)} dx \ge \sqrt[3]{36V^2\pi}.
$$
Under the given condition, {\em equality} can be reached on one but never on
both sides of the above inequality at the same time. It is clear that
$$
\mbox{(III.b)}\hspace{3cm} S\ge\sqrt[3]{36V^2\pi}.
$$
By this it seems to me that the following theorem is proven rigorously:

{\it The ball has a smaller surface than any other body of same volume which
surface is formed by a finite number of finite pieces characterized as an
algebraic surface in every point.}
\vspace{.1cm}

\centerline{\S 5}

\vspace{.1cm}
The previous examination includes any finite body ${\mathfrak U}$ which surface
is formed by a finite number of analytical surfaces.

The condition "the surface of ${\mathfrak U}$ is formed from a finite number
of {\it analytical} surfaces" is sufficient but, as one can easily be convinced,
not necessary for the conclusion that "the analyzed body ${\mathfrak U}$ has,
if it is not a ball itself, a greater surface than a ball of same volume". The
theorem mentioned in the beginning of this manuscript holds at all times even
without this limiting condition, if the surface of ${\mathfrak U}$ is formed
by a finite number of finite pieces each one having at every point a {\em
unique tangential plane} which changes continuously with the point location.

To prove  this, we just have to show that for every body ${\mathfrak U}$ that
was constructed this way and is not a ball, there exists a {\em polyhedron}
${\mathfrak U}^*$ bounded by a finite number of planar surfaces which has the
{\em same} volume but a {\em smaller} surface area than ${\mathfrak U}$.

The conclusion that the body ${\mathfrak U}$ has a greater surface area than a
ball with the same volume is proven by applying the previously performed
analysis (\S\S 1-4) on the body ${\mathfrak U}^*$ a fortiori.

Now, there exists an abstract in Steiner's paper "Ueber Maximum und Minimum bei
den Figuren in der Ebene, auf der Kugelfl\"{a}che und im Raume \"uberhaupt"
[Gesammelte Werke Band II. Seite 300-306] that provides the means to construct
for every non-spherical body ${\mathfrak U}$ with the mentioned properties, the
surface $S$ and the volume $V$, another body ${\mathfrak U}'$ which volume
$V'$ is {\em equal} to $V$ whereas the surface area $S'$ is {\em smaller} than
$S$.

If we consider such a body ${\mathfrak U}'$ and introduce two non-zero variables
$\epsilon$ and $\eta$, not subjected to any restrictions due to their smallness,
there is an infinite number of ways to construct a {\em polyhedron}
${\mathfrak U}''$ bounded by planar surfaces. We construct ${\mathfrak U}''$
in such a way that if we denote its surface area by $S''$ and the volume by
$V''$  then the differences $S'' - S'$ and $V''-V'$ in absolute values are
smaller than $\epsilon S'$ and $\eta V'$, respectively. If we now construct a
polyhedron ${\mathfrak U}''$ that is {\em similar} to ${\mathfrak U}^*$, which
volume $V^*$ is {\em equal} to the volume of body ${\mathfrak U}$, then the
surface area $S^*$ of this polyhedron is smaller than
$$
(1+\epsilon)\sqrt[3]{(1-\eta)^{-2}}\;S'.
$$
If we choose $\epsilon$ and $\eta$ from the beginning such that
$$
\frac{(1+\eta)^3}{(1-\eta)^2} < \left(\frac{S}{S'}\right)^3,
$$
then we find that $S^* < S$.

We conclude that it is possible, as asserted, for any non-spherical body
${\mathfrak U}$ with the mentioned properties, to construct a polyhedron
${\mathfrak U}^*$ bounded by a finite number of plane surfaces, which has the
{\em same} volume and a {\em smaller} surface area than ${\mathfrak U}$.

We conclude that out of all bodies of the same volume with surfaces of the
mentioned properties, the ball has the smallest surface area.
\end{document}